\documentclass[a4paper]{amsart}

\usepackage[a4paper, margin=2cm]{geometry}
\usepackage[british]{babel}

\usepackage{amsmath}
\usepackage{amssymb}
\usepackage{calc}
\usepackage{enumitem}
\usepackage{graphicx}
\usepackage{url}
\usepackage{xcolor}
\usepackage{mathtools}
\usepackage{hyperref}

\hypersetup{
  colorlinks   = true, 
  urlcolor     = {blue!50!black}, 
  linkcolor    = {blue!50!black}, 
  citecolor   = {red!50!black} 
}

\numberwithin{equation}{section}

\usepackage{cleveref}

\theoremstyle{plain}
\newtheorem{thm}{Theorem}[section]

\newtheorem{prop}[thm]{Proposition}
\newtheorem{lem}[thm]{Lemma}

\newtheorem*{lem*}{Lemma}
\newtheorem*{cor*}{Corollary}
\newtheorem*{thm*}{Theorem}
\newtheorem*{fact*}{Fact}

\theoremstyle{definition}
\newtheorem{defn}[thm]{Definition}
\newtheorem{eg}[thm]{Example}
\newtheorem{notn}[thm]{Notation}
\newtheorem{qn}[thm]{Question}

\newtheorem*{notn*}{Notation}
\newtheorem*{term*}{Terminology}
\newtheorem*{qn*}{Question}

\theoremstyle{remark}
\newtheorem{rem}[thm]{Remark}

\DeclareMathOperator{\Aut}{Aut}

\DeclareMathOperator{\Forb}{Forb}

\DeclareMathOperator{\Nminus}{N_{_{\!--}}}

\DeclareMathOperator{\sign}{sign}

\DeclareMathOperator{\sw}{sw}

\newcommand{\N}{\mathbb{N}}
\newcommand{\Q}{\mathbb{Q}}
\newcommand{\R}{\mathbb{R}}

\newcommand{\Z}{\mathbb{Z}}

\newcommand{\mb}[1]{\mathbb{#1}}
\newcommand{\mc}[1]{\mathcal{#1}}

\newcommand{\sub}{\subseteq}

\newcommand{\tld}[1]{\tilde{#1}}

\newcommand{\ra}{\rightarrow}

\newcommand{\Mod}[1]{\ (\mathrm{mod}\ #1)}

\newcommand{\eps}{\varepsilon}
\newcommand{\emp}{\varnothing}
\newcommand{\Sym}{\mathrm{S}}
\newcommand{\Alt}{\mathrm{A}}
\newcommand{\op}{\mathrm{op}}

\def\ind{\mathrel{\raise0.2ex\hbox{\ooalign{\hidewidth$\vert$\hidewidth\cr\raise-0.9ex\hbox{$\smile$}}}}}

\newcommand{\Fr}{Fra\"{i}ss\'{e} }
\newcommand{\Frn}{Fra\"{i}ss\'{e}}

\counterwithout{equation}{section}

\allowdisplaybreaks

\author{Rob Sullivan}
\address{\parbox{\linewidth}{Rob Sullivan\\
Independent Researcher\\
East London,\\
United Kingdom
}}
\email{robertsullivan1990+maths@gmail.com}

\thanks{The author was funded by Project 24-12591M of the Czech Science Foundation (GA\v{C}R).}



\title{Ultrahomogeneous hypertournaments}

\date{\today}

\begin{document}

\begin{abstract}
    A $k$-hypertournament is a higher-arity generalisation of a tournament, consisting of a set together with an irreflexive $k$-ary relation such that the automorphism group of each substructure of size $k$ is the alternating group $\mathrm{A}_k$. We show that, for each $k \geq 3$, there are continuum-many non-isomorphic countable ultrahomogeneous $k$-hypertournaments, answering a question of Cherlin.
\end{abstract}

\maketitle

\section{Introduction}

We consider the following first-order relational structures, which are higher-arity generalisations of tournaments:

\begin{defn} \label{d:Tk}
    Let $k \geq 2$. A \emph{$k$-hypertournament} $T$ consists of a set together with an irreflexive $k$-ary relation $R^T$, such that the automorphism group of each $k$-element substructure of $T$ is the alternating group $\Alt_k$.
    
    For example, a $2$-hypertournament is a tournament in the usual sense (an oriented complete graph), and a $3$-hypertournament is a structure where each triple of points has a cyclic orientation anticlockwise or clockwise.
\end{defn}

Hypertournaments under this definition have been considered in \cite{Che19}, \cite{Che22b}, \cite{CHKN21}, \cite{dlNGS25}, \cite{KSW27}, \cite{BSWY26}, \cite{Est26}, and particular examples of hypertournaments were earlier considered in \cite{Ben96}, \cite{BC00}, without being named as such. (Note that there is another definition of hypertournament in the combinatorial literature dating back to Assous (\cite{Ass86}), where the $k$-ary relation is defined to hold for exactly one enumeration of each $k$-set; we do not use this definition.)

Recall that a first-order structure $M$ is \emph{ultrahomogeneous} if every isomorphism between finitely generated substructures extends to an automorphism of $M$. We assume the reader is familiar with the basic theory of these -- see \cite{Hod93}, \cite{Mac11} for more background. We refer to a countably infinite ultrahomogeneous structure as a \emph{\Fr structure}, and we refer to the class of finitely generated structures which embed into a \Fr structure (its \emph{age}) as a \emph{\Fr class}.

The purpose of this short paper is to show the following result:

\begin{thm} \label{t: main}
    For each $k \geq 3$, there are continuum-many non-isomorphic countably infinite ultrahomogeneous $k$-hypertournaments.
\end{thm}

This answers Problem 24 in Appendix B of Cherlin's monograph \cite{Che22b}. In Subsection \ref{ss: motiv}, we state this problem (as Question \ref{q: cherlin qn}) and give motivation for it, as well as giving additional reasons for interest in Theorem \ref{t: main}. We also provide further background, giving previous examples of $k$-hypertournaments found in the literature.

\subsection{Motivation} \label{ss: motiv}

\subsubsection{Classification of ternary \Fr structures}

\begin{defn}
    Let $k \geq 2$. Let $\mc{F}$ be a set of finite $k$-hypertournaments. We write $\Forb(\mc{F})$ for the class of finite $k$-hypertournaments $T$ such that no element of $\mc{F}$ embeds in $T$, and call $\Forb(\mc{F})$ the class of \emph{finite $\mc{F}$-free $k$-hypertournaments}. If $\mc{F} = \{F\}$ for some $k$-hypertournament $F$, we call $\Forb(\mc{F})$ the class of finite $F$-free $k$-hypertournaments.
    
    If there exists $r \in \N$ such that each element of $\mc{F}$ is of size at most $r$, we say that $\Forb(\mc{F})$ is \emph{$r$-constrained}. If $\Forb(\mc{F})$ is $r$-constrained and additionally a \Fr class, we also refer to its \Fr limit as being \emph{$r$-constrained}.
\end{defn}

Cherlin asked the following:

\begin{qn}[Problem 24, Appendix B, {\cite{Che22b}} -- rephrased] \label{q: cherlin qn}
    Are there infinitely many (countably infinite, non-isomorphic) ultrahomogeneous $k$-hypertournaments for some $k$ (notably, for $k = 3$)? For $k \geq 3$, are there any that are not $(k+1)$-constrained?
\end{qn}

Note that Theorem \ref{t: main} completely answers Question \ref{q: cherlin qn}: for each $k \geq 3$, as there are only finitely many $(k+1)$-constrained ultrahomogeneous $k$-hypertournaments up to isomorphism, Theorem \ref{t: main} gives continuum-many countably infinite non-isomorphic ultrahomogeneous $k$-hypertournaments which are not $(k+1)$-constrained.

The original motivation for Question \ref{q: cherlin qn}, as stated in \cite[Subsection B.1.4]{Che22b}, was to begin a classification of \Fr structures in ternary languages; an initial step towards classification of the \Fr structures which model a particular theory is to determine how many there are. There is a long history of classification for binary languages: some kinds of \Fr structures that have been completely classified include countably infinite ultrahomogeneous posets (\cite{Sch79}), tournaments (\cite{Lac84}), graphs (\cite{LW80}), bipartite graphs with labelled parts (\cite{GGK96}), oriented graphs (\cite{Che98}) and linearly ordered graphs (\cite{Che22a}). Note that the existence of continuum-many \Fr structures modelling a certain theory does not preclude classification: Henson showed in \cite{Hen72} that there are continuum-many non-isomorphic ultrahomogeneous oriented graphs (also shown independently by Peretyat'kin in \cite{Per73}), yet a full classification of countable ultrahomogeneous oriented graphs was nonetheless provided by Cherlin in \cite{Che98}.

When considering \Fr structures in ternary languages, the simplest possible ternary language consists of a ternary relation symbol $R$; assuming the relation to be irreflexive, it is then natural to consider the ternary analogues of tournaments and graphs, namely $3$-hypertournaments and $3$-hypergraphs, which are what one respectively obtains when one imposes an additional antisymmetry or symmetry requirement on the relation. Ultrahomogeneous tournaments and graphs were among the first \Fr structures to be completely classified, giving additional motivation for beginning with their ternary analogues when approaching the classification problem for ternary \Fr structures. Some initial observations regarding countably infinite ultrahomogeneous $3$-hypergraphs can be found in \cite{AL95}. Ultrahomogeneous $3$-hypertournaments were first considered systematically by Cherlin in the unpublished notes \cite{Che19}, where the classification of the ultrahomogeneous $4$-constrained $3$-hypertournaments is given (see Example \ref{ex: 4-constr}), and Question \ref{q: cherlin qn} is first posed.

\subsubsection{Motivation: structural Ramsey theory; independence relations} \hfill

We assume in this subsection that the reader is familiar with the context and terminology of structural Ramsey theory: see, for example, \cite{HK26}. A central open question in this area (see \cite{NVT13}) is the following:

\begin{qn}
    Does every \Fr structure in a finite relational language have a coprecompact Ramsey expansion?
\end{qn}

The classification by Cherlin of the ultrahomogeneous $4$-constrained $3$-hypertournaments includes a particular structure called the \emph{generic $H_4$-free $3$-hypertournament}, which we denote by $\mb{W}$. (See Definition \ref{d: H_{k+1}} for the definition of $H_4$.) It is currently not known if $\mb{W}$ has a coprecompact Ramsey expansion. Cherlin, Hubi\v{c}ka, Kone\v{c}n\'{y} and Ne\v{s}et\v{r}il state in \cite{CHKN21} that the current standard techniques for producing Ramsey expansions from the key paper \cite{HN19} seem not to work for $\mb{W}$.

The continuum-many ultrahomogeneous $k$-hypertournaments that we produce in this paper for each $k \geq 3$ are in a certain sense a generalisation of $\mb{W}$: we define a notion of a \emph{pair-cancelling} hypertournament (Definition \ref{d: p-c}), inspired by a property satisfied by $H_4$, and each of the continuum-many $k$-hypertournaments that we construct results by taking the \Fr limit of $\Forb(\mc{P})$, where $\mc{P}$ is a collection of finite pair-cancelling $k$-hypertournaments. We do not know if any of these continuum-many $k$-hypertournaments have a coprecompact Ramsey expansion.

Another motivation for interest in Theorem \ref{t: main} comes from the investigation of independence relations on \Fr structures. A certain notion of independence relation called a \emph{stationary weak independence relation (SWIR)} was defined by Li in \cite{Li19}, \cite{Li21} as a generalisation of the notion of a stationary independence relation (SIR) due to Tent and Ziegler (\cite{TZ13a}). SWIRs (and in particular SIRs) have been used in determining the normal subgroups of automorphism groups of \Fr structures. The recent paper \cite{BSWY26} (which also contains a review of the area)  developed a technique where, given an $\omega$-categorical \Fr structure $M$, one determines the normal subgroups of $\Aut(M)$ as follows. First one finds an appropriate expansion $N$ of $M$ with a SWIR; secondly one uses this to determine the normal subgroups of $\Aut(N)$. Finally, this information can then be used to determine the normal subgroups of $\Aut(M)$. The development of this technique led the author to the following question:

\begin{qn} \label{q: SWIR exp}
    Does every $\omega$-categorical \Fr structure have an $\omega$-categorical \Fr expansion with a SWIR?
\end{qn}

Kaplan and Simon also asked this question in \cite[Conjecture 7.1]{KS19} with a weaker notion of independence relation called a canonical independence relation (CIR). We do not know if $\mb{W}$ or the continuum-many examples of ultrahomogeneous $k$-hypertournaments constructed in this paper have an $\omega$-categorical expansion with a SWIR or a CIR, and we do not know how to determine the normal subgroups of their automorphism groups.

See also \cite{Mig25} for further model-theoretic discussion of $\mb{W}$, where it is shown that the theory of $\mb{W}$ is strictly $\mathrm{NSOP}_4$. We do not know where the theories of the $k$-hypertournaments constructed in this paper lie with respect to the known model-theoretic dividing lines.

\subsection{Background: examples of \texorpdfstring{$k$}{k}-hypertournaments} \hfill

\medskip

We now survey prior examples of $k$-hypertournaments occurring in the literature relevant to Question \ref{q: cherlin qn}.

\begin{eg} \label{e: circular orders}
    \emph{Circular orders}. Let $\prec$ be a linear order on a set $D$. Define a ternary relation $R_\prec$ on $D$ by specifying that, for all distinct $a, b, c \in D$, we have $R_\prec(a, b, c)$ exactly when some cyclic permutation of the tuple $(a, b, c)$ is $\prec$-increasing. We call $R_\prec$ the \emph{circular order induced by $\prec$}; clearly $(D, R_\prec)$ is a $3$-hypertournament. We say that a $3$-hypertournament is a \emph{circular order} when it is a circular order induced by some linear order, and writing $Z_4$ for the unique (up to isomorphism) circular order on a $4$-element set, it is straightforward to show that the circular orders are exactly those $3$-hypertournaments $T$ where each $4$-element substructure is isomorphic to $Z_4$. We call the \Fr limit of the class of finite circular orders the \emph{generic circular order}; this is the circular order induced by the standard linear order on $\Q$. (See also \cite{Cam76}.) The concept of circular orders is relatively natural, and as such appears in many places in the literature; an early reference is \cite{Hun16}.
\end{eg}

\begin{eg} \label{ex: even 3-hts}
    \emph{Even $3$-hypertournaments}. We say that a $3$-hypertournament $T$ is \emph{even} if, for each $U \sub T$ with $|U| = 4$, for each linear order $\prec$ on $U$, the number of $\prec$-increasing edges of $T$ in $U$ is even. We define \emph{odd} $3$-hypertournaments analogously. (Given a linear order $\prec$ on $U$, one can obtain any other linear order on $U$ from $\prec$ by swapping consecutive elements; each swap preserves the parity of the number of increasing $3$-tuples, so it suffices to check the parity condition for a single linear order.) For example, circular orders are even. Babai and Cameron studied even $3$-hypertournaments in \cite{BC00}, under the name of ``oriented two-graphs". There is a one-to-one correspondence between switching classes of tournaments and even $3$-hypertournaments (see Remark \ref{r: switching}). One can use this correspondence to show that the class of finite even $3$-hypertournaments is a \Fr class; we call the \Fr limit the \emph{generic even $3$-hypertournament}. Bennett showed in \cite{Ben96} that the generic even $3$-hypertournament is a reduct of the random tournament. Cherlin showed in the unpublished notes \cite{Che19} that the only countably infinite ultrahomogeneous even $3$-hypertournaments are the generic circular order and the generic even $3$-hypertournament: given an even $3$-hypertournament $T$, taking $a \in T$, define a tournament $T_a$ on $T \setminus \{a\}$ by specifying for all distinct $v, w \in T_a$ that $v \ra w$ if $R^T(a, v, w)$; it is not difficult to show that $T$ is ultrahomogeneous if and only if $T_a$ is ultrahomogeneous and $\Aut(T) \curvearrowright T$ is transitive, and one then uses the classification of countable ultrahomogeneous tournaments due to Lachlan (\cite{Lac84}). 
\end{eg}

\begin{eg} \label{ex: 4-constr}
    \emph{Ultrahomogeneous $4$-constrained $3$-hypertournaments}. These were classified by Cherlin in the unpublished notes \cite{Che19}. Cherlin first shows that there are three isomorphism types of $3$-hypertournament on a set of four vertices. Taking this set to be $\{0, 1, 2, 3\}$, these isomorphism types have the following representatives.
    \begin{itemize}
        \item $H_4$: this has edge set equal to the $\Alt_4$-orbit of $(0, 1, 2)$, and is the unique (up to isomorphism) ultrahomogeneous $3$-hypertournament on $4$ vertices (see Definition \ref{d: H_{k+1}}).
        \item $O_4$: the \emph{odd} $3$-hypertournament on $4$ vertices. This has edge set given by taking the edge set of $H_4$ and reversing the orientation of the edge on $\{0, 1, 2\}$.
        \item $Z_4$: the circular order on $4$ vertices induced by the natural linear order $0 < 1 < 2 < 3$. (Note that this is denoted $C_4$ in \cite{CHKN21}.)
    \end{itemize}
    (One sees that these are the only isomorphism types by observing that each $3$-hypertournament on $4$ vertices can be obtained from $H_4$ by reversing some edges; as $\Aut(H_4) = \Alt_4$, for $n \leq 4$ any two $3$-hypertournaments obtained from $H_4$ by reversing $n$ edges are isomorphic, and as $H_4$, $O_4$ and $Z_4$ are isomorphic to their respective complements, it suffices to consider $n \leq 2$.)
    
    Cherlin then uses this to classify the countably infinite ultrahomogeneous $4$-constrained $3$-hypertournaments: these are exactly the generic $3$-hypertournament, the generic circular order, the generic even $3$-hypertournament and the generic $H_4$-free $3$-hypertournament. (Proof sketch: let $M$ be such a structure. If $M$ is even, then $M$ is the generic circular order or the generic even hypertournament (see Example \ref{ex: even 3-hts}). Otherwise $M$ contains a copy of $O_4$. By Ramsey's theorem, any infinite $3$-hypertournament contains an infinite circular order, so $M$ contains a copy of $Z_4$. If $M$ contains a copy of $H_4$, then $M$ is the generic $3$-hypertournament; otherwise $M$ is the generic $H_4$-free $3$-hypertournament.)
\end{eg}

We now give other examples of $k$-hypertournaments, which result from alternating multilinear forms; these are not relevant to Question \ref{q: cherlin qn}, but as it seems that $k$-hypertournaments have not been systematically surveyed before in the literature, we take the opportunity to provide brief descriptions here.

\begin{eg} Let $k \geq 2$.
    \begin{itemize}
        \item \emph{Chiral $k$-hypertournaments}. Let $D \sub \R^k$ be such that any $k$ points in $D$ are linearly independent. Define a $k$-ary relation $R$ on $D$ for distinct $v_0, \cdots, v_{k-1}$ by:
        \[R(v_0, \cdots, v_{k-1}) \text{ if } \det(v_0 \cdots v_{k-1}) > 0.\]
        As the determinant is an alternating multilinear form, the structure $(D, R)$ is a $k$-hypertournament. This example easily generalises: rather than the determinant, one may consider any alternating map to an ordered field.
        \item \emph{CC-systems}. These are closely related to chiral $3$-hypertournaments. Knuth defines a CC-system to consist of a set together with a ternary relation satisfying a certain list of axioms, intended to generalise the concept of the anticlockwise/clockwise orientation of non-collinear triples of points in $\R^2$. The first three axioms specify that a CC-system is a $3$-hypertournament -- we refer the reader to \cite{Knu92} for further details.
        \item \emph{Paley-Vandermonde $k$-hypertournaments}. Let $q$ be a power of a prime, with $q \equiv 3 \Mod{4}$. Then $-1$ is not a square in the finite field $\mb{F}_q$. Define a $k$-ary relation $R$ on $\mb{F}_q$ for distinct $v_0, \cdots, v_{k-1}$ by:
        \[R(v_0, \cdots, v_{k-1}) \text{ if } \prod_{i < j < k} (v_j - v_i) \text{ is a square.}\]
        Given distinct $v_0, \cdots, v_{k-1} \in \mb{F}_q$ and $\sigma \in \Sym_k$, we have $\prod_{i < j < k} (v_{\sigma(j)} - v_{\sigma(i)}) = \sign(\sigma) \prod_{i < j < k} (v_j - v_i)$, so $(\mb{F}_q, R)$ is a $k$-hypertournament.
    \end{itemize}
\end{eg}

\subsection{Structure of the paper and proof outline}\hfill

\medskip

In Section \ref{s: p-c}, we introduce the notion of a \emph{pair-cancelling} $k$-hypertournament. Note that $H_4$ has the following property: for any distinct $v, w \in H_4$, the remaining two points ``disagree on $vw$": that is, they can be labelled as $u, u'$ such that $R(u, v, w) \wedge \neg R(u', v, w)$. The notion of a pair-cancelling $k$-hypertournament generalises this property. We then show that, for $\mc{P}$ a collection of finite pair-cancelling $k$-hypertournaments, the class of $\mc{P}$-free $k$-hypertournaments has strong amalgamation: this is Proposition \ref{p: forb pair-cancelling}.

The remainder of the paper is devoted to constructing, for each $k \geq 3$, a countably infinite collection $\mc{A}$ of finite pair-cancelling $k$-hypertournaments which forms an antichain with respect to embeddings. Once one has such an antichain $\mc{A}$, by taking the \Fr limit of the class of $\mc{P}$-free $k$-hypertournaments for each $\mc{P} \sub \mc{A}$, one obtains $2^{\aleph_0}$ non-isomorphic ultrahomogeneous $k$-hypertournaments; this is standard, and the straightforward argument can be found in the proof of Theorem \ref{t: main} at the end of Section \ref{s: ac}.

To produce $\mc{A}$, we start by showing how to induce $k$-hypertournaments from tournaments in Section \ref{s: inducing hts}. In Section \ref{s: ac}, we take a certain infinite antichain of tournaments, and show that inducing a collection of $k$-hypertournaments from this antichain of tournaments via the technique of Section \ref{s: inducing hts} gives an infinite antichain of pair-cancelling $k$-hypertournaments, as required.

\medskip

\textbf{LLM disclosure\footnote{Commercial LLMs were used in this project in a spirit of experimentation; the author is rather unsure regarding their long-term impact on the mathematical ecosystem.}.}

LLMs were used to assist with proof discovery for this paper. At the 2026 North American ASL meeting, James Hanson, Nick Ramsey and Caroline Terry prompted ChatGPT 5.6 Sol (via Caroline Terry's mobile phone) to try and solve the open problems in Appendix B of Cherlin's monograph \cite{Che22b}. Amongst other results, Sol produced an initial proof of a restricted form of the main Theorem \ref{t: main} in the case $k = 3$, solving the $k=3$ case of Problem 24 of \cite{Che22b}. Caroline Terry then sent the output produced by Sol regarding this list of problems to Gabe Conant, who likewise used Sol to extract the solution to the $k=3$ case of Problem 24 and to clean up the broad structure of the proof.

The proof given by Sol for $k = 3$ involved two key ideas: the definition of a pair-cancelling $3$-hypertournament, and a seemingly ex nihilo construction of an antichain of pair-cancelling $3$-hypertournaments. Sol did not cite any sources. Gabe Conant mentioned the result to Alberto Miguel G\'{o}mez, who noted that the definition of pair-cancelling given by Sol was a generalisation of the definition of $H_4$.

Gabe Conant sent Sol's proof for $k = 3$ to the present author, after Alberto Miguel G\'{o}mez mentioned the result to him. The author observed that the pair-cancelling $3$-hypertournaments produced by Sol were just those sum-induced from the antichain of tournaments given by Henson (see Definition \ref{d: T_n}). This led the author to conjecture a proof strategy for Theorem \ref{t: main} for general $k \geq 3$: generalise the definition of pair-cancelling, and use the $k$-hypertournaments sum-induced from the Henson tournament antichain to produce an antichain of pair-cancelling $k$-hypertournaments. The author then used Claude Fable 5 to assist with some of the technical details of this proof strategy in Lemmas \ref{l: pc}, \ref{l: dir identities} and \ref{l: from k to 3} - this involved a lengthy ``conversation" between Fable and the author. Fable suggested a complicated proof of Proposition \ref{p: antichain}, but the author realised that the proof given by Sol for $k = 3$ could be used instead, resulting in the proofs of Lemmas \ref{l: Q H_4} and \ref{l: less than m(m-5)}. After obtaining a rough proof of the main result, the author wrote the paper from scratch without the use of LLMs.

\textbf{Acknowledgements.}

The author would like to thank Gabe Conant, James Hanson, Nick Ramsey and Caroline Terry for their generosity in providing him with the initial proof of the $k = 3$ case of Theorem \ref{t: main} produced by ChatGPT 5.6 Sol, as well as Alberto Miguel G\'{o}mez for mentioning this result to him in the first place.

\section{Pair-cancelling hypertournaments} \label{s: p-c}

\begin{defn} \label{d: p-c}
    Let $T$ be a $k$-hypertournament. We say that $T$ is \emph{pair-cancelling} if $|T| \geq k+1$ and, for all distinct $v, w \in T$, for each linear order $\prec$ on $T \setminus \{v, w\}$, there are $\prec$-increasing $(k-2)$-tuples $\bar{u}, \bar{u}'$ in $T \setminus \{v, w\}$ with $R(\bar{u}, v, w) \wedge \neg R(\bar{u}', v, w)$.
\end{defn}

\begin{defn} \label{d: H_{k+1}}
    We define a $k$-hypertournament $H_{k+1}$ on $k+1$ vertices $c_0, \cdots, c_k$ as follows: for each $k$-tuple $\bar{v}$ of distinct vertices, letting $w$ be the vertex of $H_{k+1}$ not contained in $\bar{v}$, we define $R(\bar{v})$ exactly when $(\bar{v}, w)$ is an even permutation of $(c_0, \cdots, c_k)$. (Equivalently: $R^{H_{k+1}}$ is the $\Alt_{k+1}$-orbit of $(c_0, \cdots, c_{k-1})$.)

    We immediately have $\Aut(H_{k+1}) \cong \Alt_{k+1}$. Also, it is easy to see that $H_{k+1}$ is ultrahomogeneous: any bijection of $(k-1)$-sets extends to a $k$-hypertournament isomorphism of substructures of size $k$, which in turn extends to an element of $\Aut(H_{k+1})$. In fact $H_{k+1}$ is the unique (up to isomorphism) ultrahomogeneous $k$-hypertournament on $k+1$ vertices: given an ultrahomogeneous $k$-hypertournament $T$ on $c_0, \cdots, c_k$, by ultrahomogeneity $\Aut(T)$ contains all $3$-cycles, so $\Aut(T) \cong \Alt_{k+1}$, and as $\Alt_{k+1}$ acts transitively on $k$-sets we have $T \cong H_{k+1}$.
\end{defn}

\begin{eg} \label{e: Hk+1 pair cancelling}
    The $k$-hypertournament $H_{k+1}$ is pair-cancelling, as we now show. Given $v, w \in H_{k+1}$ and a linear order $\prec$ on $H_{k+1} \setminus \{v, w\}$, let $\bar{u} = (u_0, \cdots, u_{k-2})$ be a $\prec$-increasing enumeration of $H_{k+1} \setminus \{v, w\}$, and for each $i < k-1$ let $\bar{t}^{(i)}$ be the $(k-2)$-tuple obtained from $\bar{u}$ by omitting $u_i$ in the $\prec$-increasing enumeration. Then $R(\bar{t}^{(i)}, v, w)$ exactly when $(\bar{t}^{(i)}, v, w, u_i)$ is an even permutation of $(c_0, \cdots, c_k)$, by the definition of $H_{k+1}$. We may obtain $(\bar{t}^{(0)}, v, w, u_0)$ from $(\bar{t}^{(1)}, v, w, u_1)$ by swapping $u_0, u_1$, and so the $(k-2)$-tuples $\bar{t}^{(0)}$, $\bar{t}^{(1)}$ are as required.
\end{eg}

\begin{eg}
    Let $n \in \N$. Let $D = \{0, \cdots, n-1\}$. Define a $k$-hypertournament $\Lambda$ on $D$ by specifying $R(v_{\sigma(0)}, \cdots, v_{\sigma(k-1)})$ for all $<$-increasing $\bar{v} \in (D)^k$ and for all $\sigma \in \Alt_k$ (where $<$ denotes the standard linear order on $D$). It is straightforward to check that $\Lambda$ does not contain any pair-cancelling $k$-hypertournament.
\end{eg}

\begin{prop} \label{p: forb pair-cancelling}
    Let $\mc{P}$ be a collection of finite pair-cancelling $k$-hypertournaments. Then $\Forb(\mc{P})$ is a \Fr class with strong amalgamation.
\end{prop}
\begin{proof}
    Let $\mc{K} = \Forb(\mc{P})$. We show that $\mc{K}$ has strong amalgamation. By a straightforward induction, it suffices to consider the case of $A, B, C \in \mc{K}$ with $B \cap C = A$ and $B \setminus A = \{b\}$, $C \setminus A = \{c\}$. Let $\prec$ be a linear order on $A$. We define a hypertournament $D$ on $A \cup \{b,c\}$ by: $D$ extends $B, C$, and for each $\prec$-increasing $(k-2)$-tuple $\bar{a}$ in $A$ we define $R^D(\bar{a}, b, c)$, closing under even permutations. It is immediate that each pair-cancelling $k$-hypertournament contained in $D$ must lie in $B$ or in $C$, so $D \in \mc{K}$.
\end{proof}

\section{Inducing hypertournaments from tournaments} \label{s: inducing hts}

\begin{notn}
    Given a set $X$ and $k \geq 1$, we write $(X)^k$ for the subset of $X^k$ consisting of $k$-tuples with distinct coordinates.
\end{notn}

\begin{defn} \label{d:sign of ht edge}
    Let $k \geq 2$. Let $T$ be a $k$-hypertournament. We define the \emph{direction map} $\delta^T : (T)^k \to \Z/2\Z$ as follows:
    \[
        \delta^T(\bar{v}) =
        \begin{cases}
            0 & \bar{v} \in R^T\\
            1 & \bar{v} \notin R^T.
        \end{cases}
    \]
    When $T$ is clear from context we just write $\delta$ instead of $\delta^T$. It is immediate that, given two $k$-hypertournaments $T, T'$ on the same domain $D$ with $\delta^T(\bar{v}) = \delta^{T'}(\bar{v})$ for all $\bar{v} \in (D)^k$, we have $T = T'$.
\end{defn}

\begin{notn}
    We write $\Sym_k$ and $\Alt_k$ for the symmetric group and the alternating group on $\{0, \cdots, k-1\}$, and for $\sigma \in \Sym_k$ we write $\Pi(\sigma) \in \Z/2\Z$ for the parity of $\sigma$ in the usual permutation group sense:
    \[\Pi(\sigma) = |\{(m, n) \in \{0, \cdots, k-1\}^2 \mid (m < n) \wedge (\sigma(m) > \sigma(n))\}| \Mod{2}.\]
    Recall that a pair $(m, n)$ with $m < n$ and $\sigma(m) > \sigma(n)$ is called an \emph{inversion} of $\sigma$, and recall that, for each decomposition of $\sigma$ into transpositions, we have that $\Pi(\sigma)$ is equal to the parity of the number of transpositions in the decomposition.
\end{notn}

The proof of the following Lemma \ref{l: def ht from fn} is immediate.

\begin{lem} \label{l: def ht from fn}
    Let $D$ be a set, and let $k \geq 2$. Let $f : (D)^k \to \Z/2\Z$ be a function such that, for all $\sigma \in \Sym_k$ and $(v_0, \cdots, v_{k-1}) \in (D)^k$, we have
    \begin{equation} \label{e: dir map cond}
        f(v_{\sigma(0)}, \cdots, v_{\sigma(k-1)}) = f(v_{0 \vphantom{\sigma(k-1)}}, \cdots, v_{k-1 \vphantom{\sigma(k-1)}}) + \Pi(\sigma).    
    \end{equation}
    Then $T := (D, \{\bar{v} \in (D)^k \mid f(\bar{v}) = 0\})$ is a $k$-hypertournament with direction map $f$. We refer to $T$ as the $k$-hypertournament \emph{induced by} $f$.
\end{lem}

\begin{defn} \label{d: induce k-ht}
    Let $T$ be a tournament, with domain $D$. Let $k \geq 2$. Define a function $f : (D)^k \to \Z/2\Z$ by:
    \begin{equation} \label{e: sum tournament}
        f(\bar{v}) := \sum_{i < j < k} \delta^T(v_i, v_j).
    \end{equation}
    Then for all $\sigma \in \Sym_k$ and $\bar{v} \in (D)^k$ we have
    \[f(v_{\sigma(0)}, \cdots, v_{\sigma(k-1)}) = \sum_{i < j < k} \delta^T(v_{\sigma(i)}, v_{\sigma(j)}) = f(v_{0 \vphantom{\sigma(k-1)}}, \cdots, v_{k-1 \vphantom{\sigma(k-1)}}) + \Pi(\sigma),\]
    so $f$ satisfies (\ref{e: dir map cond}) in Lemma \ref{l: def ht from fn}, and thus induces a $k$-hypertournament with direction map $f$, which we denote by $\Sigma_k(T)$. We call $\Sigma_k(T)$ the $k$-hypertournament \emph{sum-induced} by $T$, and write $\delta^T_k = \delta^{\Sigma_k(T)}$ for its direction map. (Note that $\Sigma_2(T) = T$.) We write $\delta^{}_k$ instead of $\delta^T_k$ when $T$ is clear from context.
\end{defn}

\subsection{Further remarks on sum-induced \texorpdfstring{$k$}{k}-hypertournaments}

We now give a number of parenthetical remarks on sum-induced $k$-hypertournaments for further context; these remarks are not used in the proof of the main Theorem \ref{t: main}, and can be skipped on a first reading of the paper.

\begin{rem} \label{r: switching}
    The author first learnt of the idea of sum-induced $k$-hypertournaments from Jeroen Winkel in the course of their work on the paper \cite{BSWY26}. In the case $k = 3$, this idea appears earlier: the equation (\ref{e: sum tournament}) appears in \cite[Section 2]{BC00} (written multiplicatively and with a sign change), and is used to show that there is a one-to-one correspondence between switching classes of tournaments and even $3$-hypertournaments. We briefly outline this in our notation for the convenience of the reader. Given a tournament $T$ and $U \sub T$, the \emph{switching} $\sw(T, U)$ of $T$ at $U$ is the tournament produced from $T$ by reversing the orientation of all edges with one vertex in $U$ and one vertex in $T \setminus U$; it is easy to check that $\Sigma_3(\sw(T, U)) = \Sigma_3(T)$. We say that two tournaments $T, T'$ on the same domain are \emph{switching-equivalent} if one is a switching of the other at some subset; it is not difficult to check that switching-equivalence is an equivalence relation, and so $T \mapsto \Sigma_3(T)$ is a well-defined map from the class of switching-classes of tournaments to the class of even $3$-hypertournaments. Given an even $3$-hypertournament $\tld{T}$ and $a \in \tld{T}$, define a tournament $\tau_a(\tld{T})$ by specifying:
    \begin{itemize}
        \item $v \ra w$ for all $v, w \in \tld{T} \setminus \{a\}$ with $R^{\tld{T}}(a, v, w)$;
        \item $a \ra v$ for all $v \in \tld{T} \setminus \{a\}$.
    \end{itemize}
     It is easy to check that, given an even $3$-hypertournament $\tld{T}$ and $a \in \tld{T}$, we have $\Sigma_3(\tau_a(\tld{T})) = \tld{T}$, and given a tournament $T$ and $a \in T$, we have $\tau_a(\Sigma_3(T)) = \sw(T, \Nminus(a))$ (where $\Nminus(a)$ denotes the in-neighbourhood of $a$), giving the required correspondence.
\end{rem}

\begin{rem} \label{r: even}
    Suppose that $k$ is odd. Generalising Example \ref{ex: even 3-hts}, we say that a $k$-hypertournament $T$ is \emph{even} if, for each $U \sub T$ with $|U| = k+1$, for each linear order $\prec$ on $U$, the number of $\prec$-increasing edges of $T$ in $U$ is even. (As in Example \ref{ex: even 3-hts}, it suffices to check the parity condition for a single linear order on $U$.)
    
    Let $T$ be a tournament. Then, for all $U \sub \Sigma_k(T)$ with $|U| = k+1$, taking a linear order $\prec$ on $U$ and writing $u_0 \prec \cdots \prec u_k$ for the $\prec$-increasing enumeration of $U$, we have
    \[\sum_{i_0 < \cdots < i_{k-1} \leq k} \delta^T_k(u_{i_0}, \cdots, u_{i_{k-1}}) = \sum_{i < j \leq k} (k-1)\delta^T(u_i, u_j) = 0.\]
    So $\Sigma_k(T)$ is an even $k$-hypertournament for all odd $k \geq 3$.
\end{rem}

\begin{rem}
    Let $T$ be an infinite tournament. Let $k > 3$. Then Lemma \ref{l: from k to 3} gives that $\Sigma_k(T)$ is not ultrahomogeneous, via the following proof. Suppose for a contradiction that $\Sigma_k(T)$ is ultrahomogeneous. As $\Sigma_3(T)$ is infinite, by the classical infinite Ramsey theorem it contains an infinite circular order, and so there are $v_0, v_1, v_2, v_3 \in \Sigma_3(T)$ such that $R^{\Sigma_3(T)}(v_i, v_{i+1}, v_{i+2})$ for all $i \in \Z/4\Z$. As $k > 3$, by ultrahomogeneity of $\Sigma_k(T)$ there is $g \in \Aut(\Sigma_k(T))$ with $g(v_3) = v_3$ such that $g$ acts on $\{v_0, v_1, v_2\}$ as the $3$-cycle $(v_0 v_1 v_2)$. By Lemma \ref{l: from k to 3} we have that $g$ is an isomorphism from $\Sigma_3(T)$ to $\Sigma_3(T)$ or $\Sigma_3(T)^\op$; the latter case cannot occur as $g$ setwise-fixes $\{v_0, v_1, v_2\}$, so $g \in \Aut(\Sigma_3(T))$. But $g$ acts on $\{v_0, v_1, v_2, v_3\}$ as $(v_0 v_1 v_2)(v_3)$, so does not preserve the circular order induced on $\{v_0, v_1, v_2, v_3\}$ by $\Sigma_3(T)$, contradiction.

    Note that in the case $k = 3$ it is possible for $\Sigma_3(T)$ to be ultrahomogeneous: the $3$-hypertournaments sum-induced by $(\Q, <)$ and the random tournament are the generic circular order and the generic even $3$-hypertournament respectively.
\end{rem}

\begin{rem}
    Definition \ref{d: induce k-ht} generalises as follows. Let $l \geq 2$, and let $T$ be an $l$-hypertournament with domain $D$. Let $k > l$. Define a function $f : (D)^k \to \Z/2\Z$ by $f(\bar{v}) = \sum_{i_0 < \cdots < i_{l-1} < k} \delta^T(v_{i_0}, \cdots, v_{i_{l-1}})$. One can then easily check that $f$ induces a $k$-hypertournament exactly when $\binom{k-2}{l-2} \equiv 1 \Mod{2}$.
\end{rem}

\section{An antichain of pair-cancelling hypertournaments} \label{s: ac}

\begin{defn} \label{d: T_n}
    Let $n \geq 1$. Define $D_n = \{0, \cdots, n-1\}$, and let $L_n$ be the tournament on $D_n$ defined by: for all $u, v \in D_n$, define $u \ra v$ if $u < v$. (Here $<$ is the standard integer order.) Let $\gamma \sub D_n^3$ be the circular order on $D_n$ induced by $<$: that is, we have $\gamma(a, b, c)$ exactly when some cyclic permutation of $(a, b, c)$ is a $<$-increasing triple.

    Define an undirected graph $C_n$ on $D_n$ with edge set $\{\{0, 1\}, \cdots, \{n-2, n-1\}, \{n-1, 0\}\}$. We write $\sim$ for the adjacency relation of $C_n$. Define a tournament $T_n$ on $D_n$ by taking $L_n$ and reversing the directed edges which are orientations of edges of $C_n$: that is, we reverse the directed edges \[0 \ra 1,\quad 1 \ra 2,\quad \cdots,\quad (n-2) \ra  (n-1),\quad 0 \ra (n-1).\]
\end{defn}

\begin{rem}
    The collection $\{T_6, T_7, \cdots\}$ forms an antichain of finite tournaments, and a slightly modified version of the definition of $T_n$ (also giving an antichain) was used by Henson in \cite{Hen72} to construct $2^{\aleph_0}$ non-isomorphic countably infinite ultrahomogeneous oriented graphs (usually now known as the \emph{Henson digraphs}) -- once one has shown that this collection of tournaments forms an antichain, the remainder of Henson's argument is essentially the same as in the proof of Theorem \ref{t: main} at the end of Section \ref{s: ac}. The definition of $T_n$ that we use, as well as the proof that $\{T_6, T_7, \cdots\}$ is an antichain, can be found in \cite[Section 1.3]{Che98}, where an exposition of Henson's construction is given.
\end{rem}

\begin{lem} \label{l: pc}
    Let $k \geq 3$ and $n \geq 4k - 4$. Then the $k$-hypertournament $\Sigma_k(T_n)$ is pair-cancelling.
\end{lem}
\begin{proof}
    First, note that $|T_n| \geq 4k - 4 \geq k + 1$. Let $v, w \in D_n$ be distinct, and let $\prec$ be a linear order on $D_n \setminus \{v, w\}$. (Throughout this proof, recall we write $<$ for the standard order of $\N$.) Let $(u_i \mid i < k - 3)$ be a $\prec$-increasing enumeration of the final $k-3$ elements of $D_n \setminus \{v, w\}$ in $\prec$. Let $B = \{u_i \mid i < k - 3\} \cup \{v, w\}$, and let $A = D_n \setminus B$. For each $a \in A$ let $\bar{u}^{(a)} = (a, u_0, \cdots, u_{k-4})$; note that $\bar{u}^{(a)}$ is $\prec$-increasing. We will find distinct $a, a' \in A$ with $\delta_k(\bar{u}^{(a)}, v, w) + \delta_k(\bar{u}^{(a')}, v, w) \neq 0$; this suffices to show that $\Sigma_k(T_n)$ is pair-cancelling. For all distinct $a, a' \in A$ we have $\delta_k(\bar{u}^{(a)}, v, w) + \delta_k(\bar{u}^{(a')}, v, w) = \sum_{b \in B} \delta(a, b) + \delta(a', b)$; thus it suffices to show that the function $A \to \Z/2\Z$, $a \mapsto \sum_{b \in B} \delta(a, b)$ is not constant.

    For all distinct $c, d \in D_n$, we have $\delta(c, d) = 1$ if and only if exactly one of the expressions $(c > d)$, $(c \sim d)$ holds. So for all $a \in A$ we have $\sum_{b \in B} \delta(a, b) = \mu(a) + \nu(a)$, where
    \[
        \mu(a) := |\{b \in B \mid b > a\}|, \quad \nu(a) := |\{b \in B \mid b \sim a\}|.
    \]
    The function $\mu$ is constant on each $<$-interval of $D_n$ contained in $A$. To show that $a \mapsto \sum_{b \in B} \delta(a, b)$ is not constant on $A$, it therefore suffices to show that there is a $<$-interval of $D_n$ contained in $A$ where $\nu$ is not constant.

    As $|B| = k-1$ and $|A| \geq 3(k - 1)$, there are $r \in D_n$ and $l \geq 3$ with $\{r + i \Mod{n} \mid i < l\} \sub A$ and $\{r-1 \Mod{n}, r+l \Mod{n}\} \cap A = \emp$: that is, there is a ``cyclic interval" of length $l \geq 3$ in $D_n$ which is contained in $A$. For each $a$ in this cyclic interval, if $a$ is an endpoint (that is, $a = r$ or $a = r + l - 1$), then $\nu(a) = 1$; otherwise $\nu(a) = 0$. As $l \geq 3$, one of the pairs $(r, r+1)$, $(r+l-2, r+l-1)$ is $<$-increasing (in the usual non-cyclic $<$-order); we have that $\nu$ is not constant on the $<$-interval consisting of this pair, so we are done.
\end{proof}

\subsection{Reduction to the case \texorpdfstring{$k=3$}{k=3}}

\begin{notn}
    Given an $n$-tuple $\bar{x}$ and $i < n$, we write $\hat{x}_i$ for the $(n-1)$-tuple obtained from $\bar{x}$ by omitting $x_i$. That is, $\hat{x}_i = (x_0, \cdots, x_{i-1}, x_{i+1}, \cdots, x_{n-1})$.
\end{notn}

\begin{lem} \label{l: dir identities}
    Let $T$ be a tournament with domain $D$. Let $k > 3$.
    \begin{enumerate}[label=(\roman*)]
        \item \label{i: dir identity odd} If $k$ is odd, then for all $\bar{v} \in (D)^3$ and $\bar{w} \in (D)^{k-2}$ disjoint from $\bar{v}$ we have:
        \[\delta^T_3(\bar{v}) = \delta^T_{k-2}(\bar{w}) + \sum_{j < k-2} \delta^T_k(\bar{v}, \hat{w}_j).\]
        \item \label{i: dir identity even} If $k$ is even, then for all $\bar{v} \in (D)^2$ and $\bar{w} \in (D)^{k-1}$ disjoint from $\bar{v}$ we have:
        \[\delta^T(\bar{v}) = \delta^T_{k-1}(\bar{w}) + \sum_{j < k-1} \delta^T_k(\bar{v}, \hat{w}_j).\]
    \end{enumerate}
    (Recall that direction maps take values in $\Z/2\Z$.)
\end{lem}
\begin{proof}
    \ref{i: dir identity odd}: We have:
    \begin{align*}
        \delta^T_{k-2}(\bar{w}) + \sum_{j < k-2} \delta^T_k(\bar{v}, \hat{w}_j) \,&=
        \Big( \sum_{i < j < k-2} \delta^T(w_i, w_j) \Big) + (k-2)\Big(\sum_{i < j < 3} \delta^T(v_i, v_j)\Big)\\ &\qquad\quad + (k-3) \Big( \sum_{\substack{i < 3,\\j < k - 2}} \delta^T(v_i, w_j) \Big) + (k-4) \Big( \sum_{i < j < k-2} \delta^T(w_i, w_j) \Big)\\
        &= \sum_{i < j < 3} \delta^T(v_i, v_j) = \delta^T_3(\bar{v}).
    \end{align*}
    The proof of \ref{i: dir identity even} is similar, and left to the reader.
\end{proof}

\begin{defn}
    Let $k \geq 2$, and let $T$ be a $k$-hypertournament with domain $D$. Define a $k$-ary relation $R^{T^\op}$ on $D$ by: $R^{T^\op} = (D)^k \setminus R^T$. It is easy to see that $(D, R^{T^\op})$ is a $k$-hypertournament, which we denote by $T^\op$. Informally, the $k$-hypertournament $T^\op$ is obtained from $T$ by, for each $k$-set, taking the opposite orientation to that of $T$. Note that, for all $\bar{v} \in (D)^k$, we have $\delta^{T^\op}(\bar{v}) = \delta^T(\bar{v}) + 1$.
\end{defn}

\begin{lem} \label{l: from k to 3}
    Let $k > 3$. Let $T, T'$ be tournaments on domains $D, D'$ where $|D| \geq k + 5$, and let $\eta : \Sigma_k(T) \to \Sigma_k(T')$ be a $k$-hypertournament embedding. Then $\eta$ is a $3$-hypertournament embedding of $\Sigma_3(T)$ into $\Sigma_3(T')$ or $\Sigma_3(T')^\op$.
\end{lem}
\begin{proof}
    Let $S = \Sigma_3(T)$, $S' = \Sigma_3(T')$. Define $\lambda : (D)^3 \to \Z/2\Z$, $\lambda(\bar{v}) := \delta^S(\bar{v}) + \delta^{S'}(\eta(\bar{v}))$. It suffices to show that $\lambda$ is constant.

    First consider the case where $k$ is odd. Given $\bar{v} \in (D)^3$ and $\bar{w} \in (D)^{k-2}$ disjoint from $\bar{v}$, by Lemma \ref{l: dir identities}\ref{i: dir identity odd} we have
    \[
        \lambda(\bar{v}) = \delta^T_{k-2}(\bar{w}) + \Big(\sum_{j < k-2} \delta^T_k(\bar{v}, \hat{w}_j)\Big) + \delta^{T'}_{k-2}(\eta(\bar{w})) + \Big(\sum_{j < k-2} \delta^{T'}_k(\eta(\bar{v}), \eta(\hat{w}_j))\Big),
    \]
    and as $\eta$ is a $k$-hypertournament embedding we have $\sum_{j < k-2} \delta^T_k(\bar{v}, \hat{w}_j) = \sum_{j < k-2} \delta^{T'}_k(\eta(\bar{v}), \eta(\hat{w}_j))$; hence
    \begin{equation} \label{e: lambda}
        \lambda(\bar{v}) = \delta^T_{k-2}(\bar{w}) + \delta^{T'}_{k-2}(\eta(\bar{w})).
    \end{equation}
    As $|D| \geq k + 5 > k + 4$, for all $\bar{u}, \bar{v} \in (D)^3$ there is $\bar{w} \in (D)^{k-2}$ disjoint from both $\bar{u}$ and $\bar{v}$, and so by (\ref{e: lambda}) we have $\lambda(\bar{u}) = \lambda(\bar{v})$. Thus $\lambda$ is constant as required.

    The case where $k$ is even is similar, using Lemma \ref{l: dir identities}\ref{i: dir identity even} and taking $\bar{w} \in (D)^{k-1}$ disjoint from $\bar{u}$ and $\bar{v}$ (here we use $|D| \geq k + 5$), and we leave this to the reader.
\end{proof}

\subsection{The case \texorpdfstring{$k=3$}{k=3}}

\begin{defn}
    Let $Q$ be the (undirected) graph on four vertices $v_0, v_1, v_2, v_3$ with edges $v_0 v_1$, $v_1 v_2$. 
    
    We denote an undirected path on $n$ vertices by $P_n$. With this notation, the graph $Q$ consists of $P_3$ together with an isolated vertex.
\end{defn}

\begin{term*}
    Below, a \emph{copy of $H_4$} refers to a substructure isomorphic to $H_4$. When we count copies of $H_4$ in a structure, we count the number of substructures isomorphic to $H_4$, rather than the number of embeddings.
\end{term*}

\begin{lem} \label{l: Q H_4}
    Let $n \geq 5$.
    \begin{enumerate}[label=(\roman*)]
        \item \label{i: H_4 iff Q} Let $V \sub D_n$, $|V| = 4$. Then $\Sigma_3(T_n)|_V \cong H_4$ if and only if $C_n|_V \cong Q$.
        \item \label{i: n(n-5) copies} There are $n(n-5)$ copies of $H_4$ in $\Sigma_3(T_n)$.
    \end{enumerate}
\end{lem}
\begin{proof}
    \ref{i: H_4 iff Q}: Define $\eps : (D_n)^2 \to \Z/2\Z$,
    \[
        \eps(u, v) =
        \begin{cases}
            1 & uv \text{ is an edge of } C_n\\
            0 & \text{otherwise}.
        \end{cases}
    \]
    For all $v_0, v_1, v_2 \in D_n$ with $\gamma(v_0, v_1, v_2)$, we have
    \begin{equation} \label{e: eps}
        \delta_3(v_0, v_1, v_2) = \sum_{i < j < 3} \eps(v_i, v_j).    
    \end{equation}
    (Recall from Definition \ref{d: T_n} that $\gamma$ is the circular order induced by $<$.) Thus the rotation $r \mapsto r + 1 \Mod{n}$ is an automorphism of $\Sigma_3(T_n)$, and so, as $n \geq 5$, by rotating if necessary we may assume $V = \{v_0, v_1, v_2, v_3\}$ with $v_0 < v_1 < v_2 < v_3$ and $v_0 \not\sim v_3$ (recall that $\sim$ denotes adjacency in $C_n$).

    It is straightforward to check that $\Sigma_3(T_n)|_V \cong H_4$ if and only if
    \begin{equation} \label{e: H4 pattern}
        (\delta_3(v_0 v_1 v_2), \delta_3(v_0 v_1 v_3), \delta_3(v_0 v_2 v_3), \delta_3(v_1 v_2 v_3)) \in \{(0, 1, 0, 1), (1, 0, 1, 0)\}.
    \end{equation}
    By (\ref{e: eps}), we have
    \[
        (\delta_3(v_0 v_1 v_2), \delta_3(v_0 v_1 v_3), \delta_3(v_0 v_2 v_3), \delta_3(v_1 v_2 v_3)) = (\eps(v_0 v_1) + \eps(v_1 v_2), \eps(v_0 v_1), \eps(v_2 v_3), \eps(v_1 v_2) + \eps(v_2 v_3)),
    \]
    and this together with (\ref{e: H4 pattern}) gives the result, via a straightforward check.

    \ref{i: n(n-5) copies}: By \ref{i: H_4 iff Q}, it suffices to show there are $n(n-5)$ copies of $Q$ in $C_n$: this follows as there are $n$ choices for the midpoint of $P_3$ in $C_n$, and then for each choice of $P_3$ there are $n-5$ possible choices of isolated vertex.
\end{proof}

\begin{lem} \label{l: less than m(m-5)}
    Let $n > m > 6$. Let $V \sub D_n$, $|V| = m$. Then the number of copies of $H_4$ in $\Sigma_3(T_n)|_V$ is strictly less than $m(m-5)$.
\end{lem}
\begin{proof}
    By Lemma \ref{l: Q H_4}\ref{i: H_4 iff Q}, it suffices to show that the number of copies of $Q$ in $C_n|_V$ is strictly less than $m(m-5)$. As $V \neq D_n$, the graph $C_n|_V$ is a free amalgam of paths. Denote these paths by $P^{(0)}, \cdots, P^{(r-1)}$. For each $i < r$, let $d_i$ be the number of vertices of the path $P^{(i)}$, and let $q(i)$ be the number of copies of $Q$ in $C_n|_V$ where the $3$-vertex path $P_3$ in the copy lies in $P^{(i)}$. We have:
    \[
        q(i) =
        \begin{cases}
            0 & d_i \leq 2\\
            m-3 & d_i = 3\\
            2(m-4) + (d_i-4)(m-5) & d_i \geq 4.
        \end{cases}
    \]
    (The first two cases above are straightforward, and the final case results from first counting when $P_3$ is at each of the two ends of $P^{(i)}$, and then counting for the remaining copies of $P_3$ where one cannot take vertices in $P^{(i)}$ adjacent to $P_3$ as the isolated vertex.)

    We have $q(i) < d_i(m-5)$ for all $i < r$ (here we use the fact that $m > 6$), and so $\sum_{i < r} q(i) < (m-5)\sum_{i < r} d_i = m(m-5)$. As $\sum_{i < r} q(i)$ is equal to the number of copies of $Q$ in $C_n|_V$, we are done.
\end{proof}

\subsection{The main theorem}

\begin{prop} \label{p: antichain}
    Let $k \geq 3$. Let $\mc{A} = \{\Sigma_k(T_n) \mid n \geq 4k-4\}$. Then $\mc{A}$ is an antichain of pair-cancelling $k$-hypertournaments.
\end{prop}
\begin{proof}
     Let $N = 4k-4$. By Lemma \ref{l: pc}, each $\Sigma_k(T_n)$ for $n \geq N$ is pair-cancelling. Suppose for a contradiction there is an embedding $\eta : \Sigma_k(T_m) \to \Sigma_k(T_n)$ for some distinct $m, n \geq N$. This implies $m < n$. As $m \geq 4k-4$ and $k \geq 3$, we have $m \geq k + 5 \geq 8$. We have that the map $\eta$ is a $3$-hypertournament embedding from $\Sigma_3(T_m)$ to $\Sigma_3(T_n)$ or $\Sigma_3(T_n)^\op$: this follows by Lemma \ref{l: from k to 3} in the case $k > 3$, and is trivial for $k = 3$. By Lemma \ref{l: Q H_4}\ref{i: n(n-5) copies} there are $m(m-5)$ copies of $H_4$ in $\Sigma_3(T_m)$. By Lemma \ref{l: less than m(m-5)}, each subset $V \sub D_n$ of size $m$ admits strictly less than $m(m-5)$ copies of $H_4$ in $\Sigma_3(T_n)|_V$, and as $H_4^\op \cong H_4$, the same is true of $\Sigma_3(T_n)^\op|_V$. This contradicts $\eta$ being a $3$-hypertournament embedding, so no such $\eta$ exists. Thus $\mc{A}$ is an antichain.
\end{proof}

\begin{rem}
    The reader may wonder if it is possible, given $k \geq 3$, to count copies of $H_{k+1}$ in each $\Sigma_k(T_n)$, rather than reducing to the case $k = 3$. In fact, this will not work: for $k \equiv 1 \Mod{4}$ and $T$ any tournament, the $k$-hypertournament $\Sigma_k(T)$ contains no copy of $H_{k+1}$, as we now show. Let $U \sub T$ with $|U| = k+1$. Recall from Remark \ref{r: even} that $\Sigma_k(T)$ is even. So, for any linear order $\prec$ on $U$, writing $u_0 \prec \cdots \prec u_k$ for the $\prec$-increasing enumeration of $U$, we have \[\sum_{i_0 < \cdots < i_{k-1} \leq k} \delta^T_k(u_{i_0}, \cdots, u_{i_{k-1}}) = 0.\]
    Recall the definition of $H_{k+1} = \{c_0, \cdots, c_k\}$ from Definition \ref{d: H_{k+1}}, and observe that, for $i \leq k$, the tuple $(c_0, \cdots, c_{i-1}, c_{i+1}, \cdots, c_k, c_i)$ is an even permutation of $(c_0, \cdots, c_k)$ exactly when $k-i$ is even. So:
    \[\sum_{i_0 < \cdots < i_{k-1} \leq k} \delta^{H_{k+1}}(c_{i_0}, \cdots, c_{i_{k-1}}) = \frac{k+1}{2} = 1,\]
    and thus $U \not\cong H_{k+1}$ as required.
\end{rem}

We now prove the main theorem of this paper.

\begin{proof}[Proof of Theorem \ref{t: main}]
    By Proposition \ref{p: antichain}, there is a countably infinite antichain $\mc{A}$ of finite pair-cancelling $k$-hypertournaments. Let $\mc{P}_0, \mc{P}_1 \sub \mc{A}$ be distinct, with $\mc{P}_0 \setminus \mc{P}_1 \neq \emp$. Let $T \in \mc{P}_0 \setminus \mc{P}_1$. As no element of $\mc{P}_1$ embeds into $T$, we have $T \in \Forb(\mc{P}_1)$. But $T \notin \Forb(\mc{P}_0)$, so $\Forb(\mc{P}_0) \neq \Forb(\mc{P}_1)$. Hence $|\{\Forb(\mc{P}) \mid \mc{P} \sub \mc{A}\}| = 2^{\aleph_0}$. For each $\mc{P} \sub \mc{A}$, by Proposition \ref{p: forb pair-cancelling} we have that $\Forb(\mc{P})$ is a \Fr class with strong amalgamation. So, by \Frn's theorem, taking the \Fr limit of each $\Forb(\mc{P})$ for $\mc{P} \sub \mc{A}$ gives $2^{\aleph_0}$ non-isomorphic ultrahomogeneous countably infinite $k$-hypertournaments.
\end{proof}

\bibliographystyle{alpha}
\bibliography{references}

\end{document}